\documentclass[12pt]{amsart}

\title{A Statistical Approach to Persistent Homology}
\author{Peter Bubenik}
\address{Cleveland State University, Department of Mathematics, 2121 Euclid Ave. RT 1515, Cleveland OH 44115-2214, USA}
\email{p.bubenik@csuohio.edu}
\thanks{This research was partially funded by the
    Swiss National Science Foundation grant 200020-105383.}
\author{Peter T. Kim}
\address{Department of Mathematics and Statistics, University of Guelph,
Guelph, Ontario N1G 2W1 Canada}
\email{pkim@uoguelph.ca}
\thanks{This research was partially funded by NSERC grant OGP46204.}
\date{\today}

\usepackage{amsmath}
\usepackage{amsthm}
\usepackage{amssymb}
\usepackage{ifpdf}

\ifpdf
\usepackage[pdftex]{graphicx}
\DeclareGraphicsExtensions{.pdf}
\usepackage{hyperref}
\else
\usepackage[dvips]{graphicx}
\DeclareGraphicsExtensions{.eps}
\usepackage[dvipdfm]{hyperref}
\fi

\newtheorem{thm}{Theorem}[section]
\newtheorem{lemma}[thm]{Lemma}
\newtheorem{prop}[thm]{Proposition}

\theoremstyle{definition}
\newtheorem{defn}[thm]{Definition}

\theoremstyle{remark}
\newtheorem{rem}[thm]{Remark}

\numberwithin{equation}{section}

\newcommand{\beq}{\begin{equation}}
\newcommand{\eeq}{\end{equation}}

\newcommand {\M} {\ensuremath {\mathcal{M}} }

\newcommand {\Z} {\ensuremath {\mathbb{Z}} }

\newcommand {\R} {\ensuremath {\mathbb{R}} }
\newcommand {\eR} {\ensuremath {\overline{\mathbb{R}}} }
\newcommand {\RP} {\ensuremath {\mathbb{RP}} }
\newcommand {\F} {\ensuremath {\mathbb{F}} }

\newcommand {\isom} {\ensuremath {\cong} }

\newcommand {\maxf} {\max(f_{\kappa})}
\newcommand {\minf} {\min(f_{\kappa})}
\newcommand {\cF} {{\mathcal F}}

\newcommand {\cC} {{\mathcal {C}}}

\newcommand {\Cd} {(C,d)}
\newcommand {\homoteq} {\approx}

\DeclareMathOperator{\Id}{Id}
\DeclareMathOperator{\im}{im}
\DeclareMathOperator{\Const}{Const}
\DeclareMathOperator{\tr}{tr}
\DeclareMathOperator{\diag}{diag}

\begin{document}

\maketitle


\begin{abstract}
  Assume that a finite set of points is randomly sampled from a
  subspace of a metric space.  Recent advances in computational
  topology have provided several approaches to recovering the
  geometric and topological properties of the underlying space.  In
  this paper we take a statistical approach to this problem. We assume
  that the data is randomly sampled from an unknown probability
  distribution.  We define two filtered complexes with which we can
  calculate the persistent homology of a probability distribution.
  Using statistical estimators for samples from certain families of
  distributions, we show that we can recover the persistent homology
  of the underlying distribution.
\end{abstract}

\section{Introduction}

There is growing interest in characterizing topological features of
data sets.  Given a finite set, sometimes called \emph{point cloud
  data (PCD)}, that is randomly sampled from a subspace $X$ of some metric
space, one hopes to recover geometric and topological properties of
$X$.  Using random samples, P. Niyogi, S. Smale and S. Weinberger
\cite{niyogiSmaleWeinberger} show how to recover the homology of
certain submanifolds. In \cite{chazalCohen-SteinerLieutier} the homotopy-type of certain compact
subsets is recovered.

A finer descriptor, developed by H. Edelsbrunner, D. Letscher, A.
Zomorodian and G. Carlsson, is that of \emph{persistent homology}
\cite{edelsbrunnerLetscherZomorodian, zomorodianCarlsson:computingPH}.  While it
is not a homotopy invariant, it is stable under small
changes~\cite{cohen-steinerEdelsbrunnerHarer}.  Using the PCD and the
metric, one can construct a filtered simplicial complex which
approximates the unknown space
$X$~\cite{deSilvaCarlsson,czcg:persistenceBarcodesForShapes}.
This leads naturally to a spectral sequence. What is unusual, is that
the homology of the start of the spectral sequence is uninteresting,
and so is what it converges to. Nevertheless, the intermediate
homology, called \emph{persistent homology} is of interest. It can be
described using \emph{barcodes}, which are analogues of the Betti
numbers.

The aim of this paper is to take a statistical approach to these
ideas. We assume that the data is sampled from a manifold with respect
to a probability distribution. Given such a distribution, we construct
two filtered chain complexes: the \emph{Morse complex}, and the
\emph{\v{C}ech complex}. For most of the distributions we consider,
these complexes are related by Alexander duality. Using persistent
homology, one can calculate the corresponding Betti barcodes, which
provide a topological description of the distribution. In the case of
the \v{C}ech complex we define a Betti--$0$ function. We apply to these
methods to several parametric families of distributions: the von
Mises, von Mises-Fisher, Watson and Bingham distributions on $S^{p-1}$
and the matrix von Mises distribution on $SO(3)$.

Given a sample, it is assumed that the underlying distribution is
unknown, but that it is one of a parametrized family. We use
statistical techniques to estimate the parameter. These are then
used to estimate the barcodes. As a result, we prove that we can
recover the persistent homology of the underlying distribution.

\begin{thm}
  Let $x_1, \ldots, x_n$ be a sample from $S^{p-1}$ according to the
  von Mises--Fisher distribution with fixed concentration parameter
  $\kappa \geq 0$. Given the sample, let $\hat{\kappa}$ be the maximum
  likelihood estimator for $\kappa$ (which is given by formula
  \eqref{est-kappa}). Let $\beta_{\kappa}$ and $\beta_{\hat{\kappa}}$
  denote the Betti barcodes for the persistent homology of the
  densities associated with $\kappa$ and $\hat{\kappa}$ using either
  the Morse or the \v{C}ech filtration. Finally let $E(\cdot)$
  denote the expectation, and $\mathcal{D}$ denote the barcode
  metric (see Definition~\ref{def:barcodeMetric}). Then,
  \begin{equation*}
    E (\mathcal{D}(\beta_{\hat{\kappa}},\beta_{\kappa})) \leq C(\kappa) n^{-1/2},
  \end{equation*}
  as $n \to \infty$, for some constant $C(\kappa)$.
\end{thm}

We also show that the classical theory of spacings \cite{pyke:spacings}
can be used to calculate the exact expectations of the Betti barcodes
for samples from the uniform distribution on $S^1$ together with their
asymptotic behavior.

As part of results, we show that the Morse filtrations of our
distributions each correspond to a relative CW-structure for the
underlying spaces. The von Mises and von Mises-Fisher distributions
correspond to the decomposition $S^{p-1} \approx * \cup_* D^{p-1}$,
the Watson distribution corresponds to $S^{p-1} \approx S^{p-2}
\cup_{\Id \amalg -\Id} (D^{p-1} \amalg D^{p-1})$, and the Bingham
distribution corresponds to $S^{p-1} \approx * \cup_{\Id \amalg -\Id}
(D^1 \amalg D^1) \cup_{\Id \amalg -\Id} (D^2 \amalg D^2) \cup \ldots
\cup_{\Id \amalg -\Id} (D^{p-1} \amalg D^{p-1})$. Finally, the Morse
filtration on the matrix von Mises distribution on $SO(3)$ corresponds
to the decomposition $\RP^2 \cup_f D^3$ where $f:S^2 \to \RP^2$
identifies antipodal points. Interestingly, the last decomposition is
obtained by using the Hopf fibration $S^0 \to S^3 \to \RP^3$.

A summary of the paper goes as follows.  In Section \ref{notation}, we
go over the background and notation used in this paper.  We review
both the statistical and the topological terminologies.  In
Section \ref{sectionPriorWork} we discuss filtrations and persistent
homology and
we develop two filtrations for densities.  In Section
\ref{sectionBettiBofS} we use the theory of spacings to give exact
estimates of the persistent homology of uniform samples on $S^1$.  In
Section \ref{sectionBarcodesOfDensities} we calculate the persistent
homology of some standard parametric families of densities on
$S^{p-1}$ and $SO(3)$.  In Section \ref{statestimation} we use maximum
likelihood estimators to recover the persistent homology of the underlying
density.

\section{Background and notation}
\label{notation}

In an attempt to make this article accessible to a broad audience, we
define some of the basic statistical and topological terms we will be
using.

\subsection{Statistics}

Given a manifold $\M$ with Radon measure $\nu$, a \emph{density} is a
function $f: \M \to [0,\infty]$ such that $f d\nu$ is a
\emph{probability distribution} on $\M$ with $\int_\M f d\nu = 1$.
A common statistical example is to take $\M = \R^p$, and $d\nu$ to
be the $p-$dimensional Lebesgue measure.  A density in this case would
be a nonnegative function that integrates to unity.  We can also take
$\M = S^{p-1}$, the $(p-1)$-dimensional unit sphere, with $d\nu$ being
the $(p-1)$-dimensional spherical measure.  In this case a density is
referred to as a \emph{directional density}.  For $\M$ a compact
connected orientable Riemannian manifold, $d\nu$ would be the measure
induced by the Riemannian structure.

In statistics, we think of a family of probability densities parametrized
accordingly
\begin{equation} \label{density_par}
\left\{ f_{\vartheta} : \vartheta \in \Theta\right\}   \ \ ,
\end{equation}
where $\vartheta$ is called a \emph{parameter} and $\Theta$ is called the
\emph{parameter space}.  The parameter space $\Theta$ can be quite general
and if it is some subset of a finite-dimensional vector space, then (\ref{density_par})
is referred to as a \emph{parametric} family of densities, otherwise it is
known as a \emph{nonparametric} family of densities.  Subsequent to this,
the corresponding statistical
problem will be referred to as either a parametric statistical procedure, or, a
nonparametric statistical procedure, depending on whether we are dealing with
a parametric, or nonparametric family of densities, respectively.

Some parametric examples are in order.  Let $\M = \R^p$ and consider
the normal family of location scale probability densities,
\begin{equation} \label{normal} f_{\mu, \sigma}(x) = (2 \pi \sigma^2)^{-p/2} \exp
\left\{ \tfrac {\|x-\mu\|^2}{2\sigma^2} \right\} \ \ , \end{equation} where
$\mu, x \in \R^p$ and $\sigma^2 \in [0,\infty)$.  Letting $\vartheta =
(\mu , \sigma^2)$, we note that this parametric problem has $\Theta =
\R^p \times [0,\infty )$ as its parameter space.

If we take $\M=S^{p-1}$, a well known example of a directional
density, and one that will be used in this paper is given by
\begin{equation} \label{vmf} f_{\mu,\kappa}(x) = c(\kappa) \exp\left\{\kappa x^t
  \mu\right\}, \end{equation} where $\mu , x \in S^{p-1}$, $\kappa \in
[0,\infty)$, $c(\kappa)$ is the normalizing constant and superscript
``$t$" denotes transpose.  The distribution arising from
$f_{\mu,\kappa}$ is called the \emph{von Mises-Fisher distribution}
where this parametric problem has $\Theta = S^{p-1} \times [0,\infty
)$ as its parameter space.

Somewhat related to the above is the situation where $\M = SO(p)$, the
space of $p \times p$ rotation matrices.  Let \begin{equation} \label{mvmf}
f_{\mu,\kappa}(x) = c(\kappa) \exp\left\{\kappa {\rm tr}\, x^t
  \mu\right\}, \end{equation} where $\mu , x \in SO(p)$, $\kappa \in [0,\infty)$
and $c(\kappa)$ is the normalizing constant.  The distribution arising
from $f_{\mu,\kappa}$ is called the \emph{matrix von Mises-Fisher
  distribution} where this parametric problem has $\Theta = SO(p)
\times [0,\infty )$ as its parameter space.

A \emph{sample} $X_1, X_2, \ldots X_N$ is a sequence of independent
and identically distributed random quantities on $\M$
drawn according to the density $f_{\vartheta}$ for some fixed but unknown
$\vartheta \in \Theta$.  The parameter of interest would be the fixed but unknown
parameter $\vartheta$, or, more generally, some transformation $\tau(\vartheta)$
thereof.  Statistically, we want to find an estimator
${\tilde \tau} = {\tilde \tau}(X_1, \ldots , X_N)$ of
$\tau(\vartheta)$. Given some metric $\gamma$ on $\tau(\Theta)$, the
performance of the estimator is evaluated relative to this metric in
expectation with respect to the joint probability density of the sample,
\begin{equation} \label{expectation}
E_{\vartheta}\gamma\left({\tilde \tau}, \tau \right)
= \int_{\M}\cdots \int_{\M}\gamma\left({\tilde \tau}, \tau\right)f_{\vartheta}
\cdots f_{\vartheta} d\nu \cdots d\nu \ \ ,
\end{equation}
where the above represents an $N-$fold integration and
$\vartheta \in \Theta$.
Thus the relative merit of one estimator over another estimator can be evaluated
using (\ref{expectation}) in a statistical decision theory context, see~\cite{berger:statisticalDecisionTheory}.

There are a wide variety of different distributions for a given
manifold, as well as sample spaces that are different manifolds.
References that discuss these topics can be found in the books by
Mardia and Jupp~\cite{mardiaJupp:book} and
Chikuse~\cite{chikuse:book}.  Furthermore, although nonparametric
statistical procedures on compact Riemannian manifolds are available, \cite{hendriks, efromovich, angersKim, kimKoo},
in this paper we will deal with parametric statistical procedures.

\subsection{Topology} \label{sectionBackgroundTopology}

Let $R$ be a commutative ring with identity. (In fact, we will only be
interested in cases where $R$ is a field, in which case $R$-modules are
vector spaces and $R$-module morphisms are linear maps of vector
spaces.)
\begin{defn}
A \emph{chain
  complex} over $R$ is a sequence of $R$-modules $\{C_i\}_{i \in \Z}$
together with $R$-module morphisms $d_i: C_i \to C_{i-1}$ called
\emph{differentials} such that $d_i \circ d_{i+1} = 0$. This condition
is often abbreviated to $d^2=0$.  The elements of $C_n$ are called
\emph{$n$-chains}.  This chain complex is denoted by $(C,d)$.
\end{defn}

\begin{defn}
  An (abstract) \emph{simplicial complex} $K$ is a set of finite,
  ordered subsets of an ordered set $\bar{K}$, such that
\begin{itemize}
\item the ordering of the subsets is compatible with the ordering of
  $\bar{K}$, and
\item if $\alpha \in K$ then any nonempty subset of $\alpha$ is also
  an element of $K$.
\end{itemize}
The elements of $K$ with $n+1$ elements are called $n$-simplices and
denoted $K_n$.
\end{defn}

\begin{defn} \label{defn:chainComplexOnK} Given a simplicial complex
  $K$, the \emph{chain complex} on $K$, denoted $(C_*(K),d)$ is
  defined as follows. Let $C_n(K)$ be the free $R$-module with basis $K_n$.
  We define the differential on $K_n$ and extend it to $C_n(K)$ by
  linearity.  For $[v_0, \ldots, v_n] \in K_n$ define
\[ d[v_0, \ldots, v_n] = \sum_i (-1)^i [v_0, \ldots, \hat{v}_i, \ldots, v_n],
\]
where $\hat{v}_i$ denotes that the element $v_i$ is omitted from the sequence.
\end{defn}

For $n\geq 0$, the \emph{standard $n$-simplex} is the $n$-dimensional
polytope in $\R^{n+1}$, denoted $\Delta^n$, whose vertices are given
by the standard basis vectors $e_0,\ldots ,e_n$. It is just
the convex hull of the standard basis vectors; that is \begin{equation}
\label{simplex}
\Delta^n = \left\{x = \sum_{i=0}^n a_i e_i \ \left| \ \forall i \ a_i
    \geq 0 \text{ and } \sum_{i=0}^n a_i = 1 \right. \right\}.
\end{equation}
There are inclusion maps
\begin{equation}
\label{inclusion}
\delta_i: \Delta^n \to \Delta^{n+1}
\end{equation}
(called the $i$-th face inclusion) are given
by $\delta_i(x_0,\ldots x_n) = (x_0,\ldots, x_{i-1}, 0, x_{i}, \ldots,
x_n)$ for $0 \leq i \leq n+1$.

\begin{defn} \label{defn:singularChainComplex}
Let $X$ be a topological space.
For $n\geq 0$, let $C_n(X)$ be the free $R$-module generated by the
set of continuous maps $\{\phi: \Delta^n \to X\}$.
For $n<0$, let $C_n(X) = 0$.
For $\phi: \Delta^n \to X$ let
\begin{equation} \label{boundarymaps}
d(\phi) = \sum_{i=0}^n (-1)^i \ \phi \circ \delta_i \ \in C_{n-1}(X).
\end{equation}
Extend this by linearity to an $R$-module morphism
$d: C_n(X) \to C_{n-1}(X)$.
One can check that $d^2=0$ so this defines a differential and $C_*(X) =
(\{C_n(X)\}_{n \in \Z}, d)$ is a chain complex,
called the \emph{singular chain complex}.
\end{defn}

\begin{defn} \label{defn:homology} Given a chain complex $(C,d)$, let
  $Z_k$ be the submodule given by $\{x \in C_k \ | \ dx = 0\}$ called
  the \emph{$k$-cycles}, and let $B_k$ be the submodule given by $\{ x
  \in C_k \ | \ \exists y \in C_{k+1} \text{ such that } dy = x\}$,
  called the \emph{$k$-boundaries}.  Since $d^2=0$, $d(dy)=0$ and thus
  $B_k \subset Z_k$.  The \emph{$k$-th homology} of $(C,d)$, denoted
  $H_k(C,d)$ is given by the $R$-module $Z_k/ B_k$.  The homologies
  $\{H_k(C,d)\}_{k \in \Z}$ form a chain complex with differential $0$
  denoted $H_*(C,d)$ and called the homology of $(C,d)$.  If $R$ is a
  principal ideal domain (for example, if $R$ is a field) and
  $H_k(C,d)$ is finitely generated, then $H_k(C,d)$ is the direct sum
  of a free group and a finite number of finite cyclic groups.  The
  \emph{$k$-th Betti number} $\beta_k(C,d)$ is the rank of the free
  group.  If $R$ is a field, then $\beta_k(C,d)$ equals the dimension
  of the vector space $H_k(C,d)$.  If $X$ is a topological space then
  $H_*(X)$ denotes the homology of the singular chain complex on $X$.
\end{defn}

\begin{defn}
  Two spaces $X$ and $Y$ are said to be homotopy equivalent (written
  $X \homoteq Y$) if there are maps $f:X \to Y$ and $g:Y \to X$ such
  that $g \circ f$ is homotopic to the identity map on $X$ and $f
  \circ g$ is homotopic to the identity map on $Y$.
\end{defn}

\begin{rem} \label{rem:contractible} If $X \homoteq Y$ then $H_*(X)
  \isom H_*(Y)$. So if $X$ is a \emph{contractible space} (that is, a
  space which is homotopy equivalent to a point), then $H_0(X) \isom
  R$ and $H_k(X) = 0$ for $k \geq 1$.
\end{rem}

\section{Filtrations and persistent homology} \label{sectionPriorWork}

From now on, we will assume that the ground ring is a field $\F$.

\subsection{Persistent homology} \label{sectionPersistentHomology}

In Definition~\ref{defn:homology} we showed how to calculate the
homology of a chain complex.  Given some additional information on the
chain complex, we will calculate homology in a more sophisticated way.
Namely, we will show how to calculate the \emph{persistent homology}
of a \emph{filtered chain complex}.  This will detect homology classes
which persist through a range of values in the filtration.

Let $\eR$ denote the totally ordered set of extended real numbers $\eR = \R \cup \{-\infty, \infty\}$. Then an increasing
\emph{$\eR$-filtration} on a chain complex $(C,d)$ is a sequence of
chain complexes $\{\cF_r(C,d)\}_{r \in \eR}$ such that $\cF_r(C,d)$ is
a subchain module of $(C,d)$ and $\cF_r(C,d) \subset \cF_{r'}(C,d)$
whenever $r \leq r' \in \eR$.  A chain complex together with a
$\eR$-filtration is called a \emph{$\eR$-filtered chain complex}.

For a filtered chain complex, the inclusions $\cF_j(C,d) \to
\cF_{j+l}(C,d)$ induce maps
\[
H_k(F_j(C,d)) \to H_k(F_{j+l}(C,d)).
\]
The image of this map is call the $l$-persistent $k$-th homology of
$\cF_j(C,d)$.

Let $Z^i_k = Z_k(\cF_i(C,d))$ and let $B_k^i = B_k(\cF_i(C,d))$.
Assume $\alpha \in Z^i_k$. Then $\alpha$ represents a homology class
$[\alpha]$ in $H_*(\cF_i(C,d))$.  Furthermore since $Z^i_k \subset
Z^{i'}_k$ for all $i'\geq i$, $\alpha$ also represents a homology
class in $H_*(\cF_{i'}(C,d))$, which we again denote $[\alpha]$.  One
possibility is that $[\alpha]\neq 0$ in $H_k(\cF_i(C,d))$ but
$[\alpha]= 0$ in $H_k(\cF_{i'}(C,d))$ for some $i'>i$.

Assume $\Cd$ is a chain complex with an $\eR$-filtration
${\cF}_r(\Cd)$ such that
\begin{equation} \label{eqnFiltrnCndn}
  \bigcup_{r \in \eR} {\cF}_r\Cd
  = \Cd \text{ and } \bigcap_{r \in \eR} {\cF}_r\Cd = 0.
\end{equation}
Equivalently, $\cF_{\infty}\Cd = \Cd$ and $\cF_{-\infty}\Cd = 0$.

\begin{lemma} \label{lemmar} Let $\Cd$ be a filtered chain complex
  satisfying \eqref{eqnFiltrnCndn}.  For any $n$-chain $\alpha \in
  \Cd$, there is some smallest $r \in \eR$ such that $\alpha \notin
  {\cF}_{r'}\Cd$ for all $r' < r$ and $\alpha \in {\cF}_{r''}\Cd$ for
  all $r'' > r$.
\end{lemma}

\begin{proof}
This follows from the definition of an $\eR$-filtration, the
assumption \eqref{eqnFiltrnCndn}, and the linear ordering of $\eR$.
\end{proof}

\begin{lemma} \label{lemmaHomologyInterval}
  For any $n$-cycle $\alpha \in Z_n$, the set of all $r\in \eR$ such
  that $0 \neq [\alpha] \in H_n({\cF}_r\Cd$ is either empty or is
  an interval.
\end{lemma}

\begin{proof}
  Let $\alpha \in Z_n$, and let $r_1$ be the corresponding value given by   Lemma~\ref{lemmar}.

  If there is some $\beta \in C_{n+1}$ such that $d\beta = \alpha$ then again   let $r_2$ be the corresponding value given by Lemma~\ref{lemmar}.  Since   $\beta \in {\cF}_j\Cd$ implies that $d\beta \in {\cF}_j\Cd$, it follows that   $r_2 \geq r_1$.  Thus $\alpha$ represents a nonzero homology class in   ${\cF}_r\Cd$ exactly when $r$ is in the (possibly empty) interval beginning   at $r_1$ and ending at $r_2$.  This interval contains $r_1$ if and only if   $\alpha \in {\cF}_{r_1}\Cd$, and it does not contain $r_2$ if and only if   $\beta \in {\cF}_{r_2}\Cd$.

If $\alpha$ is not a $k$-boundary then $\alpha$ represents a nonzero homology class in ${\cF}_r\Cd$ exactly when $r$ is in the interval $\{x \ | \ x \geq r_1\}$ or $\{x \ | \ x > r_1\}$.  beginning at $r_1$. Again this interval contains $r_1$ if and only if $\alpha \in {\cF}_{r_1}\Cd$.
\end{proof}

\begin{defn}
  For $\alpha \in Z_k$ define the \emph{persistence $k$-homology
    interval} represented by $\alpha$ to be the interval given by
  Lemma~\ref{lemmaHomologyInterval}.  Denote it by $I_{\alpha}$.
\end{defn}

\begin{defn} \label{defn:barcode} Define a \emph{Betti--$k$ barcode}
  to be a set of intervals\footnote{In
    Section~\ref{sectionPersistentHofD} we will see that using the
    \v{C}ech filtration, the Betti--$0$ barcode of manifolds will have
    uncountably many intervals, so we will define a more appropriate
    descriptor, the Betti--$0$ function. In
    Section~\ref{sectionBettiBofS} it will also be useful to convert
    finite Betti barcodes to functions so that we can analyze limiting
    and asymptotic behavior.}  $\{J_{\alpha}\}_{\alpha \in S \subset
    Z_k}$ such that
\begin{itemize}
\item $J_{\alpha}$ is a subinterval of $I_{\alpha}$, and
\item for all $r \in \eR$, $\{[\alpha] \ | \ \alpha \in S, \ r \in
  J_{\alpha}\}$ is an $\F$-basis for $H_k({\cF}_r\Cd)$.
\end{itemize}
We will sometimes use $\beta_k$ to denote a Betti--$k$ barcode.
\end{defn}

The set of barcodes has a
metric~\cite{czcg:persistenceBarcodesForShapes} defined as follows.

\begin{defn} \label{def:barcodeMetric} Given an interval $J$, let
  $\ell(J)$ denote its length. Given two intervals $J$ and
  $J'$, the \emph{symmetric difference}, $\Delta(J,J')$, between them
  is the one-dimensional measure of $J \cup J' - J \cap J'$. Given two
  barcodes $\{J_{\alpha}\}_{\alpha \in S}$ and
  $\{J'_{\alpha'}\}_{\alpha' \in S'}$, a \emph{partial matching}, $M$,
  between the two sets is a subset of $S\times S'$ where each $\alpha$
  and $\alpha'$ appears at most once. Define
  \begin{equation*}
    \mathcal{D}(\{J_{\alpha}\}_{\alpha \in S},
    \{J'_{\alpha'}\}_{\alpha' \in S'}) = \min_M
    \left( \sum_{(\alpha,\alpha') \in M}
      \Delta(J_{\alpha},J'_{\alpha'}) + \sum_{\alpha \notin M_1}
      \ell(J_{\alpha}) + \sum_{\alpha' \notin M_2} \ell(J'_{\alpha'}) \right),
  \end{equation*}
  where the minimum is taken over all partial matchings, and $M_i$ is
  the projection of $M$ to $S_i$.
  This defines a quasi-metric (since its value may be infinite). If
  desired, it can be converted into a metric.
\end{defn}

\subsection{Persistent homology from point cloud data}
\label{sectionPersistentHfPCD}
Let $(\M,\rho)$ be a manifold with a metric $\rho$.
Let $X = \{x_1, x_2, \ldots, x_n\} \subset \M$.
$X$ is called \emph{point cloud data}.
One would like to be able to obtain information on $\M$ from $X$.
If $X$ contains sufficiently many uniformly distributed points one may be
able to construct a complex from $X$ that in some sense reconstructs $\M$.

One such construction is the following $\eR$-filtered simplicial
complex called the \v{C}ech complex.  Recall that we are working over
a ground field $\F$. Let $\cC_*(X)$ be the largest simplicial complex
on the ordered vertex set $X$.  That is $\cC_0(X) = X$ and for $k\geq
1$, $\cC_k(X)$ consists of the ordered subsets of $X$ with $k+1$
elements.  Now filter this simplicial complex (along $\eR$) as
follows.  Given $r<0$, define $\cF^{\check{C}}_r(\cC_n(X))=0$ for all
$n$.  Let $B_r(x)$ denote the ball of radius $r$ centered at $x$.  For
$r \geq 0$ and $k\geq 1$, define $\cF^{\check{C}}_r(\cC_k(X))$ to be
the $\F$-vector space whose basis is the $k$-simplices $[x_{i_0},
\ldots, x_{i_k}]$ such that $\cap_{j=0}^k B_r(x_{i_j}) \neq 0$. We
remark that there are fast algorithms for computing
$\cF^{\check{C}}_r(\cC_k(X))$.\footnote{The balls of radius $r$
  centered at the points $\{x_{i_j}\}$ have nonempty intersection if
  and only if there is a ball of radius $r$ containing the points
  $\{x_{i_j}\}$. There are fast algorithms for the smallest enclosing
  ball problem\cite{fischerGaertnerKutz, gaertner:www}.}
$\cF^{\check{C}}_r(\cC_*(X))$ is called the $r$-\v{C}ech complex.  It
is the \emph{nerve} of the collection of balls $\{B_r(x_i)\}_{i=1}^n$,
and its geometric realization is homotopy equivalent to the union of
these balls.

A related construction is the Rips complex. For each $r$, the $r$-Rips
complex, $\cF^R_r(\cC_*(X))$, is the largest simplicial complex
containing $\cF^{\check{C}}_r(\cC_1(X))$. That is, $\cF^R_r(\cC_*(X))$
is the $\F$-vector space whose basis is the set of $k$-simplices
$[x_{i_0}, \ldots, x_{i_k}]$ such that $\rho(x_{i_j}, x_{i_{\ell}})
\leq r$ for all pairs $0 \leq j, \ell \leq k$.

Using either of these filtered chain complexes, one obtains a filtered
chain complex as follows.  Let $\Delta_*(\cC_*(X))$ be the chain
complex on $\cC_*(X)$.  Filter this over $\eR$ by letting
\[
\cF_r(\Delta_*(\cC_*(X))) = \Delta_*(\cF_r(\cC_*(X))) \text{, where }
\cF_r = \cF^{\check{C}}_r \text{ or } \cF^R_r.
\]
To simplify the notation, we write $\Delta_k(X) := \Delta_k(\cC_*(X))$.
We remark that these filtrations satisfy \eqref{eqnFiltrnCndn}:
\[
\bigcup_{r \in \eR} {\cF}_r(\Delta_*(X)) = \Delta_*(X) \text{ and }
\bigcap_{r \in \eR} {\cF}_r(\Delta_*(X)) = 0.
\]
Let $\alpha$ be an $n$-chain.
By Lemma~\ref{lemmar} we know that there is some $r \in \eR$ such that
 $\alpha \notin \cF_{r'}(\Delta_n(X))$ for all $r'<r$ and $\alpha \in
\cF_{r''}(\Delta_n(X))$ for all $r'' > r$.
In fact,

\begin{lemma} \label{lemmaRipsr}
Consider an $n$-chain, $\alpha = \sum_{i=1}^m \alpha_i
(x_{i_0},\ldots, x_{i_n})$. For the \v{C}ech filtration let
\[
r = \max_{i=1 \ldots m} \min \{ r_i \ | \ \exists x \text{ such that }
B_{r_i}(x) \ni x_{i_0}, \ldots x_{i_n} \},
\]
and for the Rips filtration let
\[
r = \max_{i=1\ldots m} \max_{j\neq k}
\rho(x_{i_j},x_{i_k}) \ \ .
\]
Then $\alpha \notin \cF_{r'}(\Delta_n(X))$ for all $r'<r$ and $\alpha \in
\cF_{r''}(\Delta_n(X))$ for all $r''\geq r$.
\end{lemma}

If $\alpha$ is an $n$-cycle then by Lemma~\ref{lemmaHomologyInterval}
there is a (possibly empty) persistence $n$-homology interval
corresponding to $\alpha$.
Applying Lemma~\ref{lemmaRipsr} to $\alpha$ and if there is some
$\beta \in \Delta_{k+1}(X)$ such that $d\beta = \alpha$, applying
Lemma~\ref{lemmaRipsr} to $\beta$, we get the following.

\begin{lemma} \label{lemmaRipsHomologyInterval}
  Given an $n$-cycle $\alpha$, the persistence $n$-homology interval
  associated to $\alpha$ is either empty or has the form $[r_1,r_2)$
  or $[r_1,\infty]$.
\end{lemma}



\subsection{Persistent homology of densities}
\label{sectionPersistentHofD}

Let $f_{\vartheta}$ be a probability density on a manifold $\M$ for
some $\vartheta \in \Theta$.  We will use $f_{\vartheta}$ to define two
increasing $\eR$-filtrations on $C_*(\M)$, the singular chain complex on
$\M$ (see Definition~\ref{defn:singularChainComplex}).

\subsubsection{The Morse filtration} \label{section:morse}

For $r \in \eR$, the \emph{excursion sets}
\begin{equation} \label{eqnMr}
\M_{\leq r} = \{ x \in \M \ | \ f_{\vartheta}(x) \leq r\},
\end{equation}
(used in Morse theory~\cite{milnor:morseTheory}) filter
$\M$ over $\eR$.
Hence they also provide an $\eR$-filtration of the singular chain
complex $C_*(\M)$,
\[
\cF^M_r(C_*(\M)) = C_*(\M_{\leq r}),
\]
which we call the \emph{Morse filtration}.
We remark that for all $k$,
\[
H_k(\cF^M_r C_*(\M)) = H_k(\M_{\leq r}).
\]

\subsubsection{The \v{C}ech filtration} \label{section:rips}

There is a dual increasing filtration to the Morse filtration which uses superlevel sets instead of sublevel sets. We modify this filtration slightly so that it mirrors the filtration on the \v{C}ech complex defined in Section~\ref{sectionPersistentHfPCD}, and we will call it the \emph{\v{C}ech   filtration}.  We do this since the filtrations on the \v{C}ech complex and the related Rips complex are the main filtrations used in computations of persistent homology.

Notice that in the \v{C}ech complex filtration all of the points in $X$, even distant outliers, appear when $r=0$. So the \v{C}ech filtration starts with all of the points of $M$ and the discrete topology, and then progressively connects the regions with decreasing density. 

For $r<0$ and all $k$, define $\cF^{\check{C}}_r(C_k(\M)) = 0$.
For $r\geq 0$, let $\cF^{\check{C}}_r(C_0(\M)) = C_0(\M)$.
Assume $k\geq 1$.
Let
\[
\Const_k = \{\phi:\Delta^k \to \M \ | \ \phi \text{ is constant} \}
\subset C_k(\M).
\]
For $0 \leq s \leq \infty$, let
\begin{equation} \label{eqnM1r}
\M_{\geq s} = \left\{m \in \M \ | \ f_{\vartheta}(m) \geq s \right\}.
\end{equation}
For $r \geq 0$, let
\begin{equation} \label{eqnFr}
\cF^{\check{C}}_r(C_k(\M)) = {\rm Const}_k + C_k(\M_{\geq \frac{1}{r}}).
\end{equation}
From this filtered chain complex we can calculate persistence $k$-homology intervals and Betti--$k$ barcodes just as in Section~\ref{sectionPersistentHfPCD}.

\begin{lemma} \label{lemmaHkFr}
For $k\geq 1$, \[H_k(\cF^{\check{C}}_r(C_*(\M))) \isom H_k(\M_{\geq \frac{1}{r}}) \ \ .\]
\end{lemma}


\begin{proof}
  By definition, $Z_k({\cF^{\check{C}}}_r C_*(\M)) = {\rm Const}_k +
  Z_k C_*(\M_{\geq \frac{1}{r}})$, and $B_k({\cF^{\check{C}}}_r
  C_*(\M)) = {\rm Const}_k + B_k C_*(\M_{\geq \frac{1}{r}})$.  So
\[
H_k({\cF^{\check{C}}}_r C_*(\M) \isom Z_k(C_*(\M_{\geq \frac{1}{r}})) /
B_k(C_*(\M_{\geq \frac{1}{r}})) = H_k(\M_{\geq \frac{1}{r}}).
\]
\end{proof}

Let $r \geq 0$.  Recall the notation of
Section~\ref{sectionPersistentHomology}: $Z^r_k =
Z_k(\cF^{\check{C}}_r(C_*(\M)))$ and $B^r_k =
B_k(\cF^{\check{C}}_r(C_*(\M))$.  To start, $Z^r_0 = \F[\M]$.  Then
$\cF^{\check{C}}_r(C_1(\M)) = \F[\{ \phi:\Delta^1 \to \M \ | \ \phi
\text{ is constant, or } \im{\phi} \subset \M_{\geq \frac{1}{r}}\}]$.

For two points $x,y \in M$, there is some map $\phi:\Delta^1 \to \M$
such that $\phi(0)=x$, $\phi(1)=y$ and $\im(\phi) \subset \M_{\geq
  \frac{1}{r}}$, in which case $d\phi = x-y$, if and only if $x$ and
$y$ are in the same path component of $\M_{\geq \frac{1}{r}}$.  Thus
\[
H_0(\cF^{\check{C}}_r(C_*(\M))) \isom \F [ \M / \sim ],
\]
where $x \sim y$ if and only if $x$ and $y$ are in the same path
component of $\M_{\geq \frac{1}{r}}$.

In the case where $\M_{\geq \frac{1}{r}}$ is path-connected,
$H_0(\cF^{\check{C}}_r(C_*(\M))) \isom \F [ \M / \M_{\geq \frac{1}{r}}
]$.  In particular $H_0(\cF^{\check{C}}_0(C_*(\M))) \isom
\F[\M/\M_{\geq \infty}]$.  Since $f_{\vartheta}$ is a probability
density, $\M_{\geq \infty}$ has measure $0$.  Therefore almost all $m
\in \M$ represent a distinct homology class in
$\cF^{\check{C}}_0(C_0(\M))$ and there are uncountably many
$0$-homology intervals.  As a result the Betti--$0$ barcode is not a
good descriptor.  In this section, we will describe how the
$0$-homology intervals can be used to describe a \emph{Betti--$0$
  function}, in the case where the density $f_{\vartheta}$ satisfies a
continuity condition.

More generally, as long as $\M - \M_{\geq \frac{1}{r}}$ is uncountable
and $\M_{\geq \frac{1}{r}}$ has countably many path components, then
almost all homology classes in $H_0(\cF^{\check{C}}_r(C_*(\M)))$ have a unique
representative.  In this case we use this as justification to consider
only those homology classes with a unique representative.

Assume that for all $r$, $\M - \M_{\geq \frac{1}{r}}$ is uncountable
and $\M_{\geq \frac{1}{r}}$ has countably many path components, and
that the following continuity condition holds for all $m \in \M$:
\begin{equation} \label{eqnContinuityCndn} \forall \epsilon > 0, \
  \exists \text{ injective } \phi: [0,1] \to \M \text{ s.t. } \phi(0)
  = m \text{ and } f(\phi(t)) > f(m)-\epsilon.
\end{equation}
This condition holds if $f_{\vartheta}$ is continuous.

\begin{lemma}
  Each $m \in M$ is a unique representative for $[m]$ for exactly
  those values of $r \in \left[0, \tfrac{1}{f_{\vartheta}(m)}\right)$
  or $r \in \left[0, \tfrac{1}{f_{\vartheta}(m)}\right]$.
\end{lemma}

\begin{proof}
  Let $m \in \M$.  Since $dm=0$, $m\in Z^r_0$ for $r\geq 0$.  Let $[m]
  \in H_*(\cF^{\check{C}}_r(C_*(\M)))$ denote the homology class
  represented by $m$.
  By definition $m \in \M_{\geq \frac{1}{r}}$ if and only if $r \geq
  \frac{1}{f_{\vartheta}(m)}$.  Thus $m$ is the unique representative
  for $[m]$ for $r < \frac{1}{f_{\vartheta}(m)}$.  By assumption, for
  any $\epsilon > 0$ there is a injective map $\phi: [0,1] \to \M$
  such that $\phi(0) = m$ and $f_{\vartheta}(\phi(t)) >
  f_{\vartheta}(m)-\epsilon$.  Then $\phi \in
  \cF^{\check{C}}_r(C_1(\M))$ where $r =
  \frac{1}{f_{\vartheta}(m)-\epsilon}$.  This implies that for any
  $\epsilon > 0$ there is a non-constant continuous map $\phi:
  \Delta^1 \to \M$ with $\phi(0)=m$ such that $\phi \in
  \cF^{\check{C}}_{\frac{1}{f_{\vartheta}(m)} + \epsilon}(C_1(\M))$.
  Hence $m$ is not a unique representative for $[m]$ for $r >
  \frac{1}{f_{\vartheta}(m)}$.  Therefore $m$ is a unique
  representative for $[m]$ for either $r \in
  \left[0,\frac{1}{f_{\vartheta}(m)}\right)$ or $r \in
  \left[0,\frac{1}{f_{\vartheta}(m)}\right]$.
\end{proof}

Before we formally define the Betti--$0$ function, we give the
following intuitive picture.  We draw each of our intervals
$\left[0,\frac{1}{f_{\vartheta}(m)}\right]$ or
$\left[0,\frac{1}{f_{\vartheta}(m)}\right)$ vertically starting at
$r=0$ and ending at $r=f_{\vartheta}(m)$.  Furthermore we order the
intervals from left to right according to their length.  In fact we
draw all of the intervals between $x=0$ and $x=1$, where the $x$-axis
is scaled according to the probability distribution
$f_{\vartheta}d\nu$.  The increasing curve traced by the tips of the
intervals will be called the Betti--$0$ function.

\begin{defn} \label{defn:betti0function} Formally, define the
  \emph{Betti--$0$ function} $\beta_0:(0,1] \times \Theta \to
  [0,\infty]$ as follows.\footnote{While our definition of $\beta_0$
    below \eqref{bb-0} is valid for $x=0$, we get
    $\beta_0(0,\vartheta) \equiv 0$. This does not provide any
    information, and is furthermore inappropriate in cases such as the
    von Mises distribution with $\kappa=0$ (see
    Section~\ref{sectionVonMises} below) where $\beta_0(x,\vartheta)$
    is constant and nonzero for $x>0$.}
  For $r \in [0,\infty]$, let
  \begin{equation} \label{eqngtheta} g_{\vartheta}(r) = \int_{\M_{\geq
        \frac{1}{r}}} f_{\vartheta} d\nu.
  \end{equation}
  Since $f_{\vartheta}$ is a probability density, $g_{\vartheta}$ is
  an increasing function $g_{\vartheta}: [0,\infty] \to [0,1]$ for
  each fixed ${\vartheta} \in \Theta$.  Also recall that $\M_{\geq
    \infty}$ has measure $0$ and by definition $\M_{\geq 0} = \M$.  So
  $g_{\vartheta}(0)=0$ and $g_{\vartheta}(\infty)=1$.  For $0 < x \leq
  1$, let
\begin{equation} \label{bb-0}
\beta_0(x,{\vartheta}) = \inf_{g_{\vartheta}(r) \geq x} r \ \ .
\end{equation}
If $g_{\vartheta}$ is continuous and strictly increasing,\footnote{In this case we can define
  $\beta_0(x,\vartheta)$ for $x \in [0,1]$.} then
\begin{equation} \label{bb-0c}
  \beta_0(x,{\vartheta}) = g_{\vartheta}^{-1}(x) \ \ ,
\end{equation}
for $\vartheta \in \Theta$.  That is,
$\beta_0(x,\vartheta)$ is the unique value of $r$ such that $\int_{M
  \geq \frac{1}{r}} f_{\vartheta} d\nu = x$.
\end{defn}

\subsubsection{Alexander duality}

The Morse and \v{C}ech filtration on $S^{p-1}$ are related by
Alexander duality.  Let $f$ be a density on $S^{p-1}$. Assume that $r
\in \im (f)$ and that $r < \sup (f)$.  Then $S^{p-1}_{f\leq r}$ is a
proper, nonempty subset of $S^{p-1}$.  Assume that $S^{p-1}_{f \leq
  r}$ is compact and a neighborhood retract. 

\begin{thm}[Alexander duality for the Morse and \v{C}ech filtrations on $S^{p-1}$]
Let $\tilde{H}$ denote reduced homology, let $\F$ be a field, and let $s=\frac{1}{r}$.
\[
\tilde{H}_i(S^{p-1}_{f > \frac{1}{s}}; \F) \isom
\tilde{H}^{p-2-i}(S^{p-1}_{f \leq r}; \F) \isom
\tilde{H}_{p-2-1}(S^{p-1}_{f\leq r}; \F).
\]
\end{thm}

\section{Expected barcodes of PCD} \label{sectionBettiBofS}

\subsection{Betti barcodes of uniform samples on $S^1$}
\label{sectionBettiUniform}

Let $f$ be the uniform density on $S^1$.  Let $X = \{X_1, \ldots X_n\}
\subset S^1$ be a sample drawn according to $f$.  $X$ is called the
point cloud data.  In this section we consider the Betti barcodes
obtained for the persistent homology of $\cF^R_*(\Delta_*(X))$ the
Rips complex on $X$ (see Section~\ref{sectionPersistentHfPCD}). The
metric we use on $S^1$ is $\frac{1}{2\pi}$ times the shortest arc length
between two points (we have normalized so that the total length of $S^1$ is one).

Before we continue, we introduce some notation.
Choose $\alpha$ such that $X_1 = e^{i \cos(\alpha)}$.
For $k = 2, \ldots n$ choose $U_k \in [0,1]$ such that
\[
X_k = e^{2\pi i (\alpha + U_k)}.
\]
We remark that each $U_k$ is uniformly distributed on $[0,1]$.  Now
reorder the $\{U_k\}$ to obtain the order statistic\footnote{Equality
  among any of the terms occurs with probability zero.}:
\[
0 < U_{n:1} < U_{n:2} < \ldots < U_{n:n-1} < 1.
\]
Let $U_{n:0} = 0$ and $U_{n:n} = 1$.
Reorder the $\{X_k\}$ as $\{X_{n:k}\}$ to correspond with the $\{U_{n:k}\}$.
Then for $1 \leq k \leq n$ define
\[
S_k = U_{n:k} - U_{n:k-1}.
\]
The set $S = \{S_1, \ldots S_n\}$ is called the set of
spacings~\cite{pyke:spacings}.
We remark that if $U_k = U_{n:j}$ with $1\leq j \leq n-1$ and take the
usual orientation of $S^1$, then the
distances from $X_k$ to its nearest backward neighbor and nearest
forward neighbor are $S_j$ and $S_{j+1}$, respectively.
Also the distance from $X_1$ to its neighbors is $S_n$ and $S_1$.
It is well known (for example, \cite{devroye}) that
\begin{lemma} \label{lemma:spacingsDistribution}
  $(S_1,\ldots,S_n)$ is uniformly distributed on the standard
  $(n-1)$-simplex $\{(x_1,\ldots x_n) | x_i\geq 0, \sum_{i=1}^n x_i =
  1\}$.
  It follows that
\[
P[S_1>a_1; \cdots; S_n>a_n] =
\begin{cases}
  (1-\sum_{i=1}^n a_i)^{n-1}& \text{if } \sum_{i=1}^n a_i < 1,\\
  0& \text{otherwise.}
\end{cases}
\]
and
\begin{equation} \label{eqnProbSimplex}
\text{(Whitworth, 1897)} \quad P(S_{n:n} > x) = \sum_{\substack{k \geq 1 \\ kx < 1}} (-1)^{k+1} (1-kx)^{n-1} \binom{n}{k}, \quad \forall x > 0.
\end{equation}
\end{lemma}

Finally, order the spacings to obtain
\[
0 < S_{n:1} < S_{n:2} < \ldots < S_{n:n-1} < 1.
\]

Now we are ready to calculate the homology in degree $0$.
Recall that $\beta_0(\cF^R_r(\Delta_*(X)))$ equals the dimension of
$H_0(\cF^R_r(\Delta_*(X))$, which equals the number of path components of
$\cF^R_r(\Delta_*(X))$.
Recall that $\cF^R_r(\Delta_0(X))$ is the empty set for $r<0$ and is the set
$X$ for $r \geq 0$.
So at $r=0$, there are (almost surely) exactly $n$ distinct homology
classes in $H_0(\cF^R_r(\Delta_*(X)))$.
Each homology class $[X_k]$ will no longer have a distinct
representative when the distance from $X_k$ to one of its neighbors is
equal to $r$.
That is each time $r$ passes one of the $S_k$ the dimension of
$H_*(\cF^R_r(\Delta_*(X)))$ decreases by one.
Therefore for $k = 0, \ldots {n-2}$,
\[
r \in \left[ S_{n:k}, S_{n:k+1} \right) \implies
\beta_0(\cF^R_r(\Delta_*(X))) = n-k.
\]
When $r \geq S_{n:n-1}$, $\cF^R_r(\Delta_*(X))$ is path connected so
$\beta_0(\cF^R_r(\Delta_*(X))) = 1$.
Translating this, we see that the Betti--$0$ barcode is the collection
of homology intervals
\[
[0, S_{n:k}) \text{ for $k = 1, \ldots {n-1}$ and $[0,\infty]$}.
\]

Finally, let us consider the homology in degree $1$.
Let
\[
\alpha =
(X_{n:1}, X_{n:2}) + \ldots + (X_{n:n-1}, X_{n:n}) + (X_{n:n},
X_{n:1}).
\]
This is a $1$-cycle in $\Delta_*(X)$.

\begin{lemma}
If $S_{n:n} \leq \frac{1}{2}$ then the Betti--$1$ barcode is the single
(possibly empty) persistence homology interval
\[
I_{\alpha} = [S_{n:n}, R), \quad \text{where } R \in [\tfrac{1}{3}, \tfrac{1}{2}),
\]
otherwise it is empty.
\end{lemma}

\begin{rem}
  If the large spacing $S_{n:n}$ is greater than or equal than
  $\frac{1}{2}$ then all of the points $X_1, \ldots X_n$ are
  concentrated on a semicircle, and $\cF^R_r(\Delta_*(X))$ does
  not contain any non-trivial $1$-cycles.  By \eqref{eqnProbSimplex},
  $P[S_{n:n} > \frac{1}{2}] = \frac{n}{2^{n-1}}$.
\end{rem}

\begin{proof}
  Assume that $S_{n:n} \leq \frac{1}{2}$.  If $r \geq S_{n:n}$,
  then $\alpha \in \cF^R_r(\Delta_1(X))$.  We claim that by
  using the definition of the Rips filtration and the geometry of
  $S^1$, $\alpha$ becomes a boundary at some $R \in [\frac{1}{3},
  \frac{1}{2}]$.  Since half the perimeter of $S^1$ is $\frac{1}{2}$, when $r\geq
  \frac{1}{2}$, $(X_i,X_j) \in \cF^R_r(\Delta_1(X))$ for all $X_i,X_j
  \in X$.  Thus when $r \geq \frac{1}{2}$ then $\cF^R_r(\Delta_*(X)) =
  \Delta_*(X)$ which is the full $(n-1)$-simplex on the
  vertices $X_1, \ldots X_n$.  In particular if $r \geq \frac{1}{2}$, then
  $\alpha$ is a boundary.

  Since $S_{n:n} < \frac{1}{2}$, the geometric realization of $\alpha$
  is a $n$-gon containing the center of $S^1$.  Thus if there is some
  $\beta = \sum \beta_{ijk}(X_i,X_j,X_k) \in
  \cF^R_r(\Delta_2(X))$ such that $d\beta=\alpha$ then for some
  $(X_i,X_j,X_k) \in \cF^R_r(\Delta_2(X))$ the geometric realization of
  $(X_i,X_j,X_k)$ contains the center of $S^1$.  The smallest $r$ for
  which this can happen is $\frac{1}{3}$.  So if $r <
  \frac{1}{3}$ then $\alpha$ cannot be a boundary.

  Thus $\alpha$ becomes a boundary when $r=R$ for some $R \in
  [\frac{1}{3},\frac{1}{2}]$.  If $S_{n:n} \geq \frac{1}{3}$ it is possible
  that $R = S_{n:n}$, and $\alpha$ is not a non-trivial boundary
  in any $\cF^R_r(\Delta_*(X)$.
\end{proof}

\begin{rem}
If $S_{n:n} < \frac{1}{3}$ then the Betti--$1$ barcode is a single
non-empty persistence homology interval.
Using \eqref{eqnProbSimplex}, $P[S_{n:n} \geq \frac{1}{3}] <
n\left(\frac{2}{3}\right)^{n-1}$.
\end{rem}

\subsection{Expected values of the Betti barcodes}

Let $U_1, \ldots U_{n-1}$ be a sample from the uniform distribution on
$[0,1]$.  Let $0 < U_{n:1} < U_{n:2} < \ldots < U_{n:n-1} < 1$ be the
corresponding order statistic.\footnote{We use $n$ here to match the
  notation of Section~\ref{sectionBettiUniform} where $\{U_1, \ldots,
  U_{n-1}\}$ is derived from $\{X_1,\ldots, X_n\} \in S^1$.}  Define
$U_{n:0} = 0$ and $U_{n:n} = 1$.  For $k = 1, \ldots n$, let $S_k =
U_{n:k} - U_{n:k-1}$.  Recall (Lemma~\ref{lemma:spacingsDistribution})
that the set of spacings $S = \{S_1, \ldots S_n\}$ is uniformly
distributed on the standard $(n-1)$-simplex.

Let $0 < S_{n:1} < \ldots < S_{n:n} < 1$ be the order statistic for
the spacings.
Then one can show~\cite[21.1.15]{shorackWellner} that

\begin{prop}
For $1\leq i \leq n$ the expected value of the spacings is given by
\[
E S_{n:i} = \frac{1}{n} \sum_{j=1}^i \frac{1}{n+1-j} = \frac{1}{n} \sum_{j=n+1-i}^n \frac{1}{j}
\]
\end{prop}

So the expected Betti--$0$ barcode is the collection of intervals
\[
\left\{ \left[0, \frac{1}{n} \sum_{j=1}^i \frac{1}{n+1-j}\right) \right\}_{i
      \in \{1, \ldots, n-1\}} \cup \{ [0, \infty] \},
\]
and the expected Betti--$1$ barcode is
\[
\left\{ \left[ \frac{1}{n} \sum_{j=1}^n \frac{1}{n+1-j}, \infty\right] \right\}.
\]

To obtain the Betti--$0$ function from the Betti--$0$ barcode let
\[
_n \tilde{\beta}_0(x,0) = E S_{n:\lceil (n-1)x \rceil}.
\]
The Betti--$0$ function is a normalized version of this $ \ _n \beta_0
(x,0) = c_n \ _n \tilde{\beta}_0 (x,0) $ so that $\int_0^1 \ _n
\beta_0 (x,0) dx = 1$. (In fact, $c_n =
  \frac{n-1}{1-ES_{n:n}}$, which for large values of $n$ is
  approximately equal to $n$.)  Thus,
\[
\ _n \beta_0 (x,0) =
\frac{c_n}{n} \sum_{j=1}^{\lceil (n-1)x \rceil} \frac{1}{n+1-j} = \frac{c_n}{n} \sum_{j=n+1-\lceil (n-1)x \rceil}^{n} \frac{1}{j}
\]

\begin{prop}
For $0<x<1$, as $n \to \infty$,
\[
\ _n \beta_0 (x,0) \to - \ln (1-x).
\]
\end{prop}

\begin{proof}
  By the definition of $c_n$, $\lim_{n\to \infty}\frac{c_n}{n} = 1$.
  The result then follows from the observation that
\[
\frac{1}{n} + \int_k^n \frac{1}{x} dx < \sum_{j=k}^n \frac{1}{j} < \frac{1}{k} + \int_k^n \frac{1}{x} dx
\]
and the fact that
\[
\lim_{n \to \infty} \ln \left( \frac{n}{n+1-\lceil (n-1)x \rceil} \right) = - \ln (1-x).
\]
\end{proof}

\begin{figure}
\begin{center}
\includegraphics[width=7cm,keepaspectratio=true]{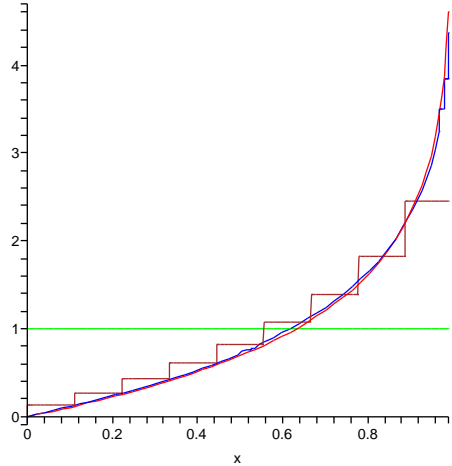}
\end{center}
\caption{Graphs of the expected Betti $0$-function for $n=10,100$ and $f(x)=-\ln(1-x)$.}
\label{figure:bettiGraph}
\end{figure}

In Figure~\ref{figure:bettiGraph}, we graph the expected Betti-$0$ functions $y= \ _{10}\beta_0(x,0)$ and $y= \ _{100}\beta_0(x,0)$ and the limiting function $y=-\ln(1-x)$. For comparison, we also graph $y=1$, the limiting function one would obtain if the spacings became relatively equal in the limit.

\section{Barcodes of certain parametric densities} \label{sectionBarcodesOfDensities}

\subsection{The von Mises distribution} \label{sectionVonMises}
Let $\M = S^1 = \{e^{i\theta} \ | \ \theta \in [-\pi,\pi)\} \subset \R^2$.
We will use this parametrization to identify $\theta \in [-\pi,\pi)$
  with an element of $S^1$.
Consider the von Mises density on $S^1$ with respect to the uniform measure,
\[ f_{\mu,\kappa}(\theta) = \tfrac{1}{I_0(\kappa)}e^{\kappa \cos(\theta -
  \mu)},  \quad \theta \in [-\pi,\pi)
\]
where $\mu \in [-\pi,\pi)$, $\kappa \in [0,\infty)$ and
$I_0(x)$ is the modified Bessel function of the first kind and
order $0$, where the general $\nu-$th order Bessel function of
the first kind is
\begin{equation}
\label{bessel}
I_{\nu}(\kappa)= \tfrac{(\kappa/2)^{\nu}}{\Gamma\left(\nu + \frac{1}{2}\right)\Gamma\left(\frac{1}{2}\right)}
\int_{-1}^1e^{\kappa t}(1-t^2)^{\nu-\frac{1}{2}}dt \ \ ,
\end{equation}
and $\Gamma (\cdot)$ denotes the gamma function.

Our homologies will be independent of $\mu$, so assume
that $\mu=0$ and so in this case the parameter $\vartheta = \kappa$.

We will filter the chain complex on $S^1$ using both the \v{C}ech and
Morse filtrations.
Recall that by \eqref{eqnM1r} and \eqref{eqnMr},
$S^1_{\geq \frac{1}{r}} = \{\theta \in S^1 \ | \ f_{\kappa}(\theta)
\geq \frac{1}{r} \}$ and
$S^1_{\leq r} = \{\theta \in S^1 \ | \ f_{\kappa}(\theta)
\leq r \}$.
Choose $\alpha_{r,\kappa} \in [-\pi,\pi)$ such that
\[
f_{\kappa}(\alpha_{r,\kappa}) = r.
\]
Specifically, let
$
\alpha_{r,\kappa} = cos^{-1}(\frac{1}{\kappa} \ln
(\frac{r}{c(\kappa)})).
$
Our calculations of the persistent homology will follow from the
following straightforward result.

\begin{lemma} \label{lemmaS1}
For $0 \leq r < \frac{1}{\max f_{\kappa}}$, $S^1_{\geq \frac{1}{r}} = \phi$,
and for  $r < \min f_{\kappa}$, $S^1_{\leq r} = \phi$.
For $\frac{1}{\max f_{\kappa}} \leq r < \frac{1}{\min f_{\kappa}}$,
\[
S^1_{\geq \frac{1}{r}} = \{ \theta \ | \ -\alpha_{\frac{1}{r},\kappa} \leq \theta \leq
\alpha_{\frac{1}{r},\kappa} \}.
\]
For $\min f_{\kappa} \leq r < \max f_{\kappa}$,
\[
S^1_{\leq r} = \{ \theta \ | \ \alpha_{r,\kappa} \leq \theta \leq 2\pi
- \alpha_{r,\kappa} \}.
\]
For $r \geq \frac{1}{\min f_{\kappa}}$, $S^1_{\geq \frac{1}{r}} = S^1$,
and for  $r \geq \max f_{\kappa}$, $S^1_{\leq r} = S^1$.
\end{lemma}

Since its analysis is simpler, we start with the Morse filtration on
$S^1$.  By Lemma~\ref{lemmaS1}, $S^1_{\leq r}$ is empty if $r < \min
f_{\kappa}$, it is contractible (see Remark~\ref{rem:contractible}) if
$\min f_{\kappa} \leq r < \max f_{\kappa}$ and it is equal to $S^1$ if
$r \geq \maxf_{\kappa}$.  It follows that the Betti--$0$ barcode for
the Morse filtration is the single interval
\[
\left[ \min f_{\kappa}, \infty \right] = \left[\tfrac{1}{ I_0(\kappa)
      e^{\kappa}}, \infty \right],
\]
the Betti--$1$ barcode is the single interval
\[
\left[ \max f_{\kappa}, \infty \right] = \left[\tfrac{e^{\kappa}}{
      I_0(\kappa)}, \infty \right],
\]
and all other Betti--$k$ barcodes are empty.

Now consider the \v{C}ech filtration on $S^1$.
We will derive a formula for the Betti--$0$ function, $\beta_0(x, \kappa)$, and
calculate the Betti--$k$ barcodes for $k>0$.

If $\kappa=0$ then $f_0 = 1$.
So for $r<1$, $S^1_{\geq \frac{1}{r}} = \emptyset$, and for $r\geq 1$, $S^1_{\geq \frac{1}{r}} = S^1$.
By definition~\eqref{eqngtheta},
\[
g_{\kappa}(r) = \begin{cases}
0 & \text{if $r < 1$},\\
1 & \text{if $r \geq 1$}.
\end{cases}
\]
So by definition \eqref{bb-0}, $\beta_0(x,0) = 1$.

For $\kappa>0$, let $\minf = \frac{1}{ I_0(\kappa)}e^{-\kappa}$ and
$\maxf=\frac{1}{ I_0(\kappa)}e^{\kappa}$.
For $r<\frac{1}{\maxf}$, $S^1_{\geq \frac{1}{r}} = \emptyset$, and for
$r\geq\frac{1}{\minf}$, $S^1_{\geq \frac{1}{r}} = S^1$.
For $\frac{1}{\maxf} \leq r < \frac{1}{\minf}$, since $f_{\kappa}$ is even and
decreasing for $\theta>0$,
\[
S^1_{\geq \frac{1}{r}} = \{ \theta \ | \ -\alpha_{r,\kappa} \leq \theta \leq \alpha_{r,\kappa}\},
\]
where $\alpha_{r,\kappa} \in (0,\pi)$ and $f_{\kappa}(\alpha_{r,\kappa}) = \frac{1}{r}$.

\begin{figure}
\begin{center}
\includegraphics[width=7cm,keepaspectratio=true]{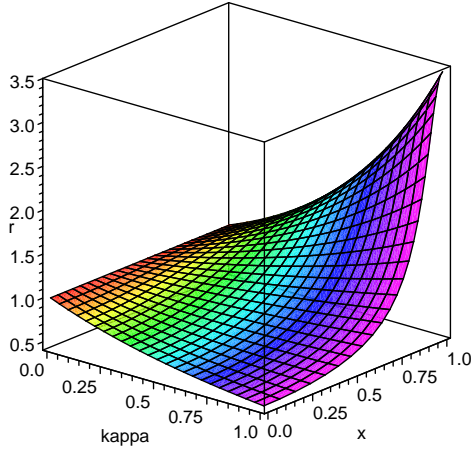}
\end{center}
\caption{Graph of the Betti $0$-function of the von Mises density for
  a range of concentration parameters}
\label{figure:vonMisesBetti0}
\end{figure}

Let $x \in [0,1]$ and assume that $\beta_0(x,\kappa)=r$.
Since $\kappa\geq 0$, $g_{\kappa}(r) = \int_{S^1_{\geq \frac{1}{r}}}f_{\kappa}(\theta)d\theta$ is continuous and strictly increasing. So,
\[
x = \int_{S^1_{\geq \frac{1}{r}}} f_{\kappa}(\theta)d\theta .
\]
Define $\alpha_{r,\kappa} \in [0,\pi]$ by the condition that $f_{\kappa}(\alpha_{r,\kappa}) =
\frac{1}{r}$.
So
\begin{equation} \label{eqn:r}
r = \frac{1}{f_{\kappa}(\alpha_{r,\kappa})}.
\end{equation}
For $\psi \in [0,\pi]$, let
\[
F_{\kappa}(\psi) = \int_0^{\psi}f_{\kappa}(\theta) d\theta.
\]
Then
\begin{equation} \label{eqn:x}
x = \int_{S^1_{\geq \frac{1}{r}}} f_{\kappa} d\nu = \int_{-\alpha_{r,\kappa}}^{\alpha_{r,\kappa}} f_{\kappa}(\theta)
d\theta = 2 F_{\kappa}(\alpha_{r,\kappa}).
\end{equation}
Since $F_{\kappa}$ is strictly increasing, it is invertible.  So
$\alpha_{r,\kappa} = F_{\kappa}^{-1}(\frac{x}{2})$.  Thus
\begin{equation} \label{betti-0} \beta_0(x,\kappa) = r =
\frac{1}{f_{\kappa}(F_{\kappa}^{-1}(\frac{x}{2}))} \end{equation} Since
$f_{\kappa}$ and $F_{\kappa}$ are smooth, by the inverse function
theorem, so is $F_{\kappa}^{-1}$.  So
\[
\beta_0(x,\kappa) = (F_{\kappa}^{-1})'\left(\frac{x}{2}\right).
\]
We remark that as $\kappa \rightarrow 0$, $\beta_0(x,\kappa)
\rightarrow 1 = \beta_0(x,0)$.
We can also describe the graph of $r=\beta_0(x,\kappa)$ parametrically
by combining \eqref{eqn:r} and \eqref{eqn:x} (see Figure~\ref{figure:vonMisesBetti0}):
\begin{equation} \label{eqn:vMh}
h_{\kappa}(t) = \left( 2 F_{\kappa}(t), \frac{1}{f_{\kappa}(t)} \right), t \in [0,\pi].
\end{equation}

For $k\geq 1$, recall that
\[
\cF^{\check{C}}_r(C_k(S^1)) = \Const_k + C_k(S^1_{\geq \frac{1}{r}}).
\]
Also recall that for $r<\frac{1}{\maxf}$, $S^1_{\geq \frac{1}{r}} = \emptyset$, for
$\frac{1}{\maxf} \leq r < \frac{1}{\minf}$, $S^1_{\geq \frac{1}{r}}$ is the arc from
$-\alpha_{r,\kappa}$ to $\alpha_{r,\kappa}$ where $f_{\kappa}(\alpha_{r,\kappa}) = \frac{1}{r}$ and for $r\geq
\frac{1}{\minf}$, $S^1_{\geq \frac{1}{r}} = S^1$.
It follows that for $k\geq 1$,
\[
H_k(\cF^{\check{C}}_r(C_*(S^1))) = \begin{cases}
\F & \text{for $k=1$ and $r\geq \frac{1}{\minf}$},\\
0 & \text{otherwise.}
\end{cases}
\]
Therefore the Betti--$1$ barcode has the single interval
\begin{equation}
\label{1-betti}
\left[\tfrac{1}{\minf},\infty\right] = \left[
  I_0(\kappa)e^{\kappa},\infty\right]
\end{equation}
and for $k>1$ the Betti--$k$ barcode is
empty.

\subsection {The von Mises-Fisher distribution}
Now consider
$\M=S^{p-1}$, $p \geq 3$ and the unimodal von Mises-Fisher density
given by
\[
f_{\mu,\kappa}(x) = c(\kappa) \exp\left\{\kappa x^t \mu\right\}, \quad x \in S^{p-1}
\]
where $\kappa \in [0,\infty)$, $\mu \in S^{p-1}$, and
\begin{equation} \label{normalizing}
  c(\kappa)=\left(\frac{\kappa}{2}\right)^{\frac{p}{2}-1}
  \frac{1}{\Gamma(\frac{p}{2}) I_{\frac{p}{2} -1}(\kappa )}
\end{equation}
is the normalizing constant with respect to the uniform measure. This
is also known as the Langevin distribution. Note that the minimum and
maximum of $f$ also do not depend on $\mu$: $\minf =
c(\kappa)e^{-\kappa}$ and $\maxf = c(\kappa)e^{\kappa}$.  In fact, by
symmetry the homologies will not depend on $\mu$.  Hence once again
take $\vartheta = \kappa$.

Consider the Morse filtration (defined in Section~\ref{section:morse})
on $S^{p-1}$.  If $r< \minf$ then $S^{p-1}_{\leq r} = \phi$ and if $r
\geq \maxf$ then $S^{p-1}_{\leq r} = S^{p-1}$.  For $\minf) \leq r <
\maxf$
\[
S^{p-1}_{\leq r} = \{ x\in S^{p-1} | x^t\mu \leq a_{r,\kappa}\},
\]
where $a_{r,\kappa} = \frac{1}{\kappa} \ln
\left(\frac{r}{c(\kappa)}\right) \in [-1,1]$.  So $S^{p-1}_{\leq r}$
is the closure of $S^{p-1}$ minus a right circular cone with vertex
$0$ and centered at $\mu$.  In particular, $S^{p-1}_{\leq r}$ is
contractible (see Remark~\ref{rem:contractible}) so
$H_0(\cF^M_r(C_*(S^{p-1}))) = \F$ and for $k\geq 1$,
$H_k(\cF^M_r(C_*(S^{p-1}))) = 0$.

Thus the Betti--$0$ barcode is the single interval $[\minf, \infty)$,
the Betti--$(p-1)$ barcode is the single interval $[\maxf, \infty)$ and all
other barcodes are empty.

Consider the \v{C}ech filtration (defined in
Section~\ref{section:rips}) on $S^{p-1}$.
For $\frac{1}{\max(f_{\kappa})} \leq r < \frac{1}{\min(f_{\kappa})}$,
\[
S^{p-1}_{\geq \frac{1}{r}} = \{x \in S^{p-1} \ | \ x^t\mu \geq a_{\frac{1}{r},\kappa}\}.
\]
So $S^{p-1}_{\geq \frac{1}{r}}$ is the intersection of $S^{p-1}$ and a right circular cone with
vertex $0$ and centered at $\mu$.
In particular for
$\frac{1}{\max(f_{\kappa})} \leq r < \frac{1}{\min(f_{\kappa})}$,
$S^{p-1}_{\geq \frac{1}{r}}$ is contractible, so for $k\geq 1$, $H_k(S^{p-1}_{\geq \frac{1}{r}})=0$.

Assume $\kappa = 0$. Then $f_0 = c(0)$, and
\[
S^{p-1}_{\geq \frac{1}{r}} = \begin{cases}
\phi & \text{if $r < \frac{1}{c(0)}$},\\
S^{p-1} & \text{if $r \geq \frac{1}{c(0)}$}.
\end{cases}
\]
Thus
\[
g_{\kappa}(r) = \begin{cases}
0 & \text{if $r < \frac{1}{c(0)}$},\\
1 & \text{if $r \geq \frac{1}{c(0)}$}.
\end{cases}
\]
Therefore $\beta_0(x,0) := \inf_{g_{\kappa}(r) \geq x} r = \frac{1}{c(0)}$.

Assume $\kappa > 0$. Then for $k=0$,
 \begin{equation} \label{eqn:g_kappa}
  x = g_{\kappa}(r) =
   \int_{S^{p-1}_{\geq \frac{1}{r}}} f_{\kappa} =
     c(\kappa) \frac{s_{p-2}}{s_{p-1}} \int_0^{\arccos
       \left(-\frac{\ln(rc(\kappa))}{\kappa}\right)} e^{\kappa \cos
       \theta} \sin^{p-2}\theta \ d\theta \ \ ,
 \end{equation}
where $s_{p-1} = \frac{2\pi^{\frac{p}{2}}}{\Gamma\left(\frac{p}{2}\right)}$.
When $\kappa > 0$, $g_{\kappa}(r)$ is continuous and strictly increasing. Hence
\begin{equation}
\label{0-betti-sphere}
\beta_0(x , \kappa ) = g_{\kappa}^{-1}(x)
\end{equation}
for $x \in [0,1]$ and $\kappa > 0$.  As we did for the von Mises
distribution~\eqref{eqn:vMh}, we can describe the graph of $r=\beta_0(x,\kappa)$ more
  explicitly using a parametric equation:
  \begin{equation} \label{eqn:vMFh}
    h_{\kappa}(t) = \left( c(\kappa) \frac{s_{p-2}}{s_{p-1}} \int_0^t e^{\kappa \cos \theta} \sin^{p-2}\theta \ d \theta, \frac{e^{-\kappa \cos t}}{c(\kappa)} \right), \quad t \in [0,\pi].
  \end{equation}

For $k\geq 1$, by Lemma~\ref{lemmaHkFr},
\[ H_k(\cF^{\check{C}}_r(C_*(S^{p-1}))) = H_k(S^{p-1}_{\geq \frac{1}{r}}) = \begin{cases}
\F & \text{ if } k=p-1 \text{ and } r \geq \frac{1}{\minf},\\
0 & \text{otherwise.}
\end{cases}
\]
Therefore for $k\geq 1$ the Betti--$k$ barcode has the single interval:
\begin{equation} \label{betti-k}
\left[\tfrac{1}{\minf},\infty\right] = \left[\tfrac{e^{\kappa}}{c(\kappa)},\infty\right]
\end{equation}
for $k=p-1$ and is empty otherwise.

\subsection{The Watson distribution} \label{sectionWatson}

Let $\M = S^{p-1}$ and consider the following bimodal distribution
\begin{equation} \label{eqnWatson}
f_{\mu,\kappa}(x) = d(\kappa) \exp \{ \kappa (x^{t}\mu)^2 \},
\end{equation}
where $\kappa\geq 0$ and $x,\mu \in S^{p-1}$,
called the \emph{Watson distribution}.
We remark that this density is rotationally symmetric, where $\mu$ is
the axis of rotation.
The minimum and maximum densities are given by
\[
\min f = d(\kappa), \quad \max f = d(\kappa) e^{\kappa}.
\]
The maximum is achieved at $x=\pm \mu$ and the minimum is achieved at
all $x$ such that $x^t \mu = 0$.

Using the Morse filtration we get the following Betti barcodes.
For $p=2$, we remark that for $r < \min f$, $S^1_{\leq r} = \phi$.
For $r = \min f$, $S^1_{\leq r}$ is two points.
As $r$ increases, these points become two arcs of increasing size,
which connect when $r = \max f$.
So the Betti--$0$ barcode consists of the two homology intervals
$[\min f, \infty] $ and $[\min f, \max f)$, and the Betti--$1$ barcode
    has the single interval $[\max f, \infty]$.
All other Betti barcodes are empty.

For $p>2$, we observe similar behavior.
When $r < \min f$, $S^{p-1}_{\leq r} = \phi$.
For $r = \min f$, $S^{p-1}_{\leq r}$ is equator which is homeomorphic
to $S^{p-2}$.
As $r$ increases, the equator expands until it reaches the poles when
$r = \max f$.
So the Betti--$0$, Betti--$(p-2)$ and Betti--$(p-1)$ barcodes each
consist of a single homology interval:
$[\min f, \infty]$, $[\min f, \max f)$, and  $[\max f, \infty]$,
      respectively.
All other Betti barcodes are empty.

Using the \v{C}ech filtration, $S^{p-1}_{\geq \frac{1}{r}}$ is either
empty, or consists of two contractible components, or is all of $S^{p-1}$.
So the Betti--$(p-1)$ barcode is the single
homology interval $[\frac{1}{\min f}, \infty]$ and the Betti--$k$
  barcodes for all other $k\geq 1$ are empty.
The Betti--$0$ function is given by $\beta_0(x,\kappa) =
g^{-1}_{\kappa}(x)$, where
\[
g_{\kappa}(r) = \int_{S^{p-1}_{\geq \frac{1}{r}}} f_{\kappa} = 2
  \frac{s_{p-2}}{s_{p-1}} \int_0^{\alpha_{\kappa}(r)} d(\kappa)
  e^{\kappa \cos^2(\theta)} \sin^{p-2}(\theta) d\theta,
\]
with $\alpha_{\kappa}(r) =
\cos^{-1}\left(\sqrt{-\frac{1}{\kappa}\ln(d(\kappa)r)}\right)$ and
$s_{p-1} = \frac{2\pi^{p/2}}{\Gamma(p/2)}$.
As with the von Mises \eqref{eqn:vMh} and von Mises-Fisher distributions \eqref{eqn:vMFh}, the Betti--$0$ function can also be described parametrically.

\subsection{The Bingham distribution}

Again let $\M = S^{p-1}$ with the probability density
\[
f_{K}(x) = d(K) \exp \{ x^t K x \}
\]
where $x \in S^{p-1} \subset \R^3$ and $K$ is a symmetric $p \times p$
matrix.
We remark that $f_{K}(x) = d(K) \exp \{ \tr K x x^t \}$.
Also, by a change of coordinates we can write $K = \diag (k_1, \ldots
k_p)$, where $k_p \geq \ldots \geq k_1$ are the eigenvalues of $K$.
Let $v_i$ be the eigenvector associated to $k_i$.

Assume that $k_p > \ldots > k_1 > 0$.
Then the minimum and maximum values of $f_K$ are given by
\[
\min f_K = d(K) e^{k_1}, \quad \max f_K = d(K) e^{k_p},
\]
and are attained at $\pm v_1$ and $\pm v_p$.

The Betti--$k$ barcodes (for $k\geq 1$) when $p=2$ are the same as for
the Watson distribution. When $p\geq 3$, the Bingham distribution
differs significantly from the Watson distribution.  For example,
the minimum of the function is attained at only $\pm v_1$ which is
certainly not homeomorphic to $S^{p-2}$.

Consider the Morse filtration.  We can calculate the Betti--$k$
barcodes inductively.  If we consider $v_p$ to be the north pole, then
there is a homotopy from $S^{p-1} - \{v_p, -v_p\}$ to $S^{p-2}$ which
collapses the sphere with missing its poles to the equator.  When $r <
k_p$, by the symmetry of $f_K$ this homotopy also gives a homotopy
from $S^{p-1}_{\leq r}$ to $S^{p-2}_{\leq r}$ where the filtration on
$S^{p-2}$ is the Morse filtration associated to the Bingham
distribution with $K = \diag (k_1, \ldots k_{p-1})$.

As a result, the Betti--$0$ barcode is given by the two homology
intervals $[d(K)e^{k_1}, \infty]$ and $[d(K)e^{k_1},
    d(K)e^{k_{2}})$.
For $1 \leq k \leq p-2$, the Betti--$i$ barcode is given by the
interval $[d(K)e^{k_{i+1}}, d(K)e^{k_{i+2}})$.
Finally, the Betti--$(p-1)$ barcode is given by the interval
$[d(K)e^{k_p}, \infty]$.

We remark that this barcode corresponds the cellular construction of
$S^{p-1}$ that repeatedly attaches northern and southern hemispheres
of increasing dimension.

For the \v{C}ech filtration we can use the same argument starting with
$v_1$. The Betti--$0$ barcode is given by the two homology intervals
$\frac{1}{d(K)} \left[e^{-k_p}, \infty\right]$ and $\frac{1}{d(K)} \left[ e^{-k_p}, e^{-k_{p-1}} \right)$. For $1 \leq i \leq
  p-2$, the Betti--$i$ barcode is given by the single interval $\frac{1}{d(K)}
  \left[e^{-k_{p-i}}, e^{-k_{p-i-1}} \right)$.
The Betti--$(p-1)$ barcode is given by the single interval $\frac{1}{d(K)} \left[ e^{-k_1}, \infty \right]$.

We remark that the correspondence between the two sets of barcodes is
a manifestation of Alexander duality.

\subsection{The matrix von Mises distribution and a Hopf fibration}
The Lie group of rotations of $\R^3$, $SO(3)$, can be given the matrix
von Mises density
\begin{equation} \label{eqnMatrixVonMises}
f_{A,\kappa}(X) = c(\kappa) \exp \left\{ \kappa \tr (X^t A) \right\},
\end{equation}
where $A \in SO(3)$ and $\kappa > 0$ is a concentration parameter.
We determine the Morse and \v{C}ech filtrations of $SO(3)$ via the
Hopf fibration $S^3 \to \RP^3$.

The special orthogonal group $SO(3)$ is diffeomorphic to the real
projective space $\RP^3$. The map $S^3 \to \RP^3$ which identifies
each point on the sphere with the one-dimensional subspace on which
it lies is a Hopf fibration whose fiber is $S^0 = \{-1,1\}$. Thus,
$S^3$ is a double-cover of $SO(3)$ (and since $S^3$ is
simply-connected, it is the universal cover).

If we represent $S^3$ with the unit quaternions and $\RP^3$ with
$SO(3)$, then the Hopf fibration above is represented by the
Cayley-Klein map $\rho: S^3 \to SO(3)$:
\[ \rho \left(
\begin{array}{c}
p_1 \\
p_2 \\
p_3 \\
p_4
\end{array} \right) = I + 2 p_1B + 2B^2 \text{, where } B = \left(
\begin{array}{ccc}
0 & -p_4 & p_3 \\
p_4 & 0 & -p_2 \\
-p_3 & p_2 & 0
\end{array} \right).
\]
We can use this map to relate the matrix von Mises
density~\eqref{eqnMatrixVonMises} on $SO(3)$ to the Watson
density~\eqref{eqnWatson} on $S^3$ by making the following
observation.
If $P = \rho(p)$ and $Q = \rho(q)$, then
\[
\tr(P^tQ) = 4 (p^t q)^2 - 1.
\]
Then if $\rho(a) = A$,
\[
\rho^{-1} \{ X \in SO(3) \ | \ f_{A,\kappa}(x) = r\} = \{ x \in S^3 \ | \
f_{a,4\kappa}(x) = kr\} 
\text{, where } k = \frac{d(4\kappa)e^{\kappa}}{c(\kappa)}.
\]
It follows that
\[
\rho^{-1}(SO(3)_{\leq r}) = S^3_{\leq kr} \text{ and }
\rho^{-1}(SO(3)_{\geq \frac{1}{r}}) = S^3_{\geq k\frac{1}{r}},
\]
where the filtration on $S^3$ is with respect to Watson density $f_{a,4\kappa}$.

Recall (Section~\ref{sectionWatson}) that for $\frac{1}{\max f} \leq
kr < \frac{1}{\min f}$, $S^3_{\geq \frac{1}{kr}}$ consists of two
contractible components. The Hopf fibration $S^3 \to \RP^3$ and
equivalently the map $\rho: S^3 \to SO(3)$ identify these two
components. So $SO(3)_{\geq \frac{1}{r}}$ is contractible. Therefore,
for the \v{C}ech filtration the Betti--$3$ barcode is the single homology
interval $[\frac{1}{\min f}, \infty)$ and all other Betti--$k$ barcodes
for $k\geq 1$ are empty. The Betti--$0$ function is identical to the
one for the Watson density on $S^3$.

For $\min f \leq kr < \max f$, $S^3_{\leq kr}$ is homotopy equivalent
(via a projection onto its equator) to $S^2$. The Hopf fibration $S^3
\to \RP^3$ restricted to the equator gives the Hopf fibration and
double cover $S^2 \to \RP^2$. The homotopy equivalences $S^3_{\leq kr}
\homoteq S^2$ induces a homotopy equivalence $SO(3)_{\leq r} \homoteq
\RP^2$. Thus for the Morse filtration, the Betti--$0$ and Betti--$3$
barcodes are the single homology intervals $[\min f, \infty)$ and
$[\max f, \infty)$ and all Betti--$k$ barcodes for $k>3$ are empty.
However, since the fundamental group and integral homology group of
degree one of $\RP^2$ are the cyclic group of order two, the Betti--$1$
and Betti--$2$ barcodes depend on the choice of the field of
coefficients $\F$. If $\F$ is a field of characteristic $0$ (e.g. the
rationals) then both are empty. However if $\F$ is the field of
characteristic two ($\Z/2\Z$), then both are the single homology
interval $[\min f, \max f)$.

\section{Statistical estimation of the Betti barcodes}
\label{statestimation}

In this section we will calculate the expected persistent homology
using statistics sampled from various densities.

\subsection{The von-Mises and von-Mises Fisher distributions}
For point cloud data $x_1, \ldots , x_n$ on $S^{p-1}$ sampled from the
von Mises-Fisher distribution (\ref{vmf}): $f_{\mu,\kappa}(x) =
c(\kappa)\exp\{\kappa x^t \mu \}$, we will give the statistical
estimators for the (unknown) parameters. We will show that these can
be used to obtain good estimates of the persistent homology of the
underlying distribution.

Letting $\bar x = \frac{1}{n} \sum_{i=1}^nx_i$ denote the sample mean,
consider the decomposition
\[ {\bar x} = \|{\bar x}\|\left(\tfrac{{\bar x}}{ \|{\bar x}\|}\right)
\ \ .
\]
The statistical estimator for $\mu$ is ${\bar x}/\|{\bar x}\|$ while
the statistical estimator for $\kappa$ is solved~\cite[Section
10.3.1]{mardiaJupp:book} by inverting $A_p(\hat \kappa) = \|{\bar
  x}\|$, where $A_p(\lambda) =
\tfrac{I_{p/2}(\lambda)}{I_{p/2-1}(\lambda)}$, and $I_\nu(\lambda)$
is the modified Bessel function of the first kind and order $\nu$.
Hence,
\begin{equation}
\label{est-kappa} {\hat \kappa} = A_p^{-1}(\|{\bar x}\|).
\end{equation}

A large sample asymptotic normality calculation for
(\ref{est-kappa}) is~\cite[Section 10.3.1]{mardiaJupp:book}
\begin{equation}
\label{asymp-mse-kappa} \sqrt{n}\left( {\hat \kappa} - \kappa
\right) \rightsquigarrow N\left(0, A_p'(\kappa)^{-1}\right),
\end{equation}
as $n \rightarrow \infty$, where $\rightsquigarrow$ means
convergence in distribution and $N(0,\sigma^2)$ stands for a
normally distributed random variable with mean 0 and variance
$\sigma^2 > 0$. 
Using this estimate of $\kappa$ we
obtain estimates for the $\beta_{\kappa}$ barcodes for the Morse and
\v{C}ech filtrations. For the Morse filtration, we estimate the
$\beta_0$ barcode and $\beta_{p-1}$ barcode to be
$[c(\hat{\kappa})e^{-\hat{\kappa}},\infty]$ and
$[c(\hat{\kappa})e^{\hat{\kappa}},\infty]$, respectively.  For the
\v{C}ech filtration, we estimate the $\beta_{p-1}$ barcode to be
$[\frac{e^{\hat{\kappa}}}{c(\hat{\kappa})},\infty]$.

Recall that the space of barcodes has a metric $\mathcal{D}$ (see
Definition~\ref{def:barcodeMetric}). Let $\beta_i^{M}(f)$ and
$\beta_i^{\check{C}}(f)$ denote the Betti--$i$ barcode for the density
$f$ using the Morse and \v{C}ech filtrations.  Then the expectations
of the distance from the estimated persistent homology to the
persistent homology of the underlying density can be bounded as follows.

\begin{thm}
  For the von Mises--Fisher distribution on $S^{p-1}$ and
  $\kappa \in [\kappa_0, \kappa_1]$, where $0 < \kappa_0 \leq
  \kappa_1 < \infty$,
  \begin{equation*}
    E (\mathcal{D}(\beta_i^M (f_{\hat{\kappa}}),\beta_i^M(f_{\kappa}))) \leq C(\kappa) n^{-1/2}
  \end{equation*}
  as $n \to \infty$ for all $i$, and
  \begin{equation*}
    E (\mathcal{D}(\beta_i^{\check{C}} (f_{\hat{\kappa}}), \beta_i^{\check{C}}(f_{\kappa}))) \leq C(\kappa) n^{-1/2}
  \end{equation*}
  as $n \to \infty$ for all $i \geq 1$, for some constant $C(\kappa)$.
\end{thm}

\begin{proof}
  Since the barcodes have a particularly simple form, we only need to
  know the barcode metric for the following case:
  \begin{equation*}
    \mathcal{D} ( \{ [a,\infty] \}, \{ [b,\infty] \} ) = |a-b|.
  \end{equation*}
  Using our previous calculations of the Betti barcodes, we have:
  \begin{eqnarray*}
    \mathcal{D} ( \beta_0^M (f_{\hat{\kappa}}), \beta_0^M(f_{\kappa}) )
    & = & |
    c(\hat{\kappa}) e^{-\hat{\kappa}} - c(\kappa) e^{-\kappa} | \\
    \mathcal{D} ( \beta_{p-1}^M (f_{\hat{\kappa}}),
    \beta_{p-1}^M(f_{\kappa}) ) & = & |
    c(\hat{\kappa}) e^{\hat{\kappa}} - c(\kappa) e^{\kappa} | \\
    \mathcal{D} ( \beta_{p-1}^{\check{C}} (f_{\hat{\kappa}}),
    \beta_{p-1}^{\check{C}}(f_{\kappa}) ) & = & |
    c(\hat{\kappa})^{-1} e^{\hat{\kappa}} - c(\kappa)^{-1} e^{\kappa} |.
  \end{eqnarray*}

  We note that the normalizing constant can be re-expressed as
  \[
  c(\kappa) = \frac{B\left(\frac{p-1}{2},\frac12\right)}{\int_{-1}^1
  e^{\kappa t}(1-t^2)^{\frac{p-3}{2}}dt}  \ \ ,
  \]
  where $B(\cdot,\cdot)$ is the beta function.
  Furthermore,
  \[
  c'(\kappa) = -B\left(\frac{p-1}{2},\frac12\right)\frac{\int_{-1}^1
  e^{\kappa t}t(1-t^2)^{\frac{p-3}{2}}dt}{\left(\int_{-1}^1
  e^{\kappa t}(1-t^2)^{\frac{p-3}{2}}dt\right)^2}  \ \
  \]
  and
  \[
  A_p'(\kappa) = 1 - A_p(\kappa)^2 - \frac{p-1}{\kappa}A_p(\kappa) \
  \ .
  \ \ \]

  For $0 \leq \kappa_0 \leq \kappa_1 < \infty$ and $\kappa \in
  \left[\kappa_0 , \kappa_1\right]$, we observe $0 < c(\kappa) ,
  |c'(\kappa)|, A_p'(\kappa) < \infty$, and by the mean value theorem,
\[E | c(\hat{\kappa}) e^{ \hat{\kappa}} - c(\kappa) e^{\kappa}|
  = E | (c(\kappa^*)+c'(\kappa^*)) e^{{\kappa}^*} ({\hat
  \kappa}-\kappa)| \ \ ,
\]
where $\kappa^*$ is a value between $\hat \kappa$ and $\kappa$.
Consequently,
\begin{eqnarray*} E | c(\hat{\kappa}) e^{
\hat{\kappa}} - c(\kappa) e^{\kappa}| &\leq& \bar{C}(\kappa)
\left\{E|{\hat \kappa}-\kappa|^2\right\}^{1/2} \\
&\leq& C(\kappa)n^{-1/2}
\end{eqnarray*}
where the first inequality is by the H\"{o}lder inequality, and the
second is by~\eqref{asymp-mse-kappa}.

Similarly,
\[E | c(\hat{\kappa}) e^{ -\hat{\kappa}} - c(\kappa) e^{-\kappa}|
  = E | (c'(\kappa^*)-c(\kappa^*)) e^{-{\kappa}^*} ({\hat
  \kappa}-\kappa)| \ \ ,
\]
and
\[E \left| \frac{e^{ \hat{\kappa}}}{c(\hat{\kappa})} - \frac{e^{\kappa}}{c(\kappa)}\right|
  = E \left| \left(\frac{c(\kappa^*)-c'(\kappa^*)}{c(\kappa^*)^2}\right) e^{{\kappa}^*} ({\hat
  \kappa}-\kappa)\right| \ \ . \qedhere
\] 
\end{proof}

Expressing the estimated $\beta_0$-function is more challenging.  For
the case of the sphere $S^2$, an exact expression can be obtained.
One can calculate that $c(\kappa) = \frac{\kappa}{\sinh(\kappa)}$, and
from \eqref{eqn:g_kappa},
\[
g_{\kappa}(r) = \frac{e^{\kappa}}{2\sinh(\kappa)} - \frac{1}{2\kappa
  r}.
\]
from which we use (\ref{0-betti-sphere}) to obtain,
\begin{equation} \label{0-betti-2}
  \beta_0(x , \kappa ) = \frac{e^{2\kappa}-1}{2\kappa[(1-x)e^{2\kappa}+x]}
\end{equation}
for $x \in (0,1]$ and $\kappa > 0$.  Notice that $\beta_0(x , \kappa )
\rightarrow 1$ as $\kappa \rightarrow 0$ and $\beta_0(x , \kappa )
\rightarrow 0$ as $\kappa \rightarrow \infty$, for all $x \in (0,1)$.
Furthermore, for (\ref{est-kappa}), \cite[9.3.9]{mardiaJupp:book}
\begin{equation}
\label{A_3}
A_3(\kappa) = \coth \kappa - \tfrac 1 {\kappa} \ \ .  \end{equation}

We have the following:
\begin{thm} For the von Mises-Fisher distribution on $S^2$, and fixed
  $\kappa > 0$,
  \[
  E \left|\left| \beta_0(x,{\hat \kappa}) -
    \beta_0(x,\kappa) \right|\right|_{\infty} \leq C(\kappa) n^{-1} \ \ ,
  \]
  as $n \rightarrow \infty$.
\end{thm}

\begin{proof}
  By the mean value theorem,
  \begin{equation} \label{eqn:mvt} \beta_0(x,{\hat \kappa})-
      \beta_0(x,\kappa) = \frac{\partial}{\partial
        \kappa} \beta_0(x,\tilde{\kappa}) ({\hat \kappa} - \kappa)
      \ ,
  \end{equation}
  where $\tilde{\kappa}$ is between $\hat{\kappa}$ and $\kappa$.
  One can calculate that
  \[
  \frac {\partial}{\partial {\kappa}}\beta_0(x,\kappa) =
  \frac{-(1-x)e^{4\kappa} + (1+2\kappa
      -2x)e^{2\kappa} + x}
  {2\kappa^2\left[(1-x)e^{2\kappa}+x\right]^2} \ \ .
  \]
  Recall that the domain of $\beta_0(x,\kappa)$ is $(0,1]$.  For $x
  \in (0,1]$, $\left| \frac{\partial}{\partial \kappa}
    \beta_0(x,\kappa) \right|$ is bounded: for instance,
  \begin{equation} \label{eqn:ddkappaBound} \left| \frac
      {\partial}{\partial {\kappa}}\beta_0(x,\kappa)\right| \leq
    \frac{e^{4\kappa} + (1+2\kappa)e^{2\kappa}+1}{2\kappa^2} \ \ .
  \end{equation}
  Combining \eqref{eqn:mvt}, \eqref{eqn:ddkappaBound}, \eqref{asymp-mse-kappa} and
  \eqref{A_3} produces the desired result.
\end{proof}

\subsection{The Watson distribution}

Recall that the Watson distribution on $S^{p-1}$ is given by
\begin{equation} \label{eqn:watson}
f_{\mu,\kappa}(x) = d(\kappa) \exp \{ \kappa (x^t \mu)^2 \} \text{, where }
\mu \in S^{p-1} \text{ and } \kappa > 0.
\end{equation}
Let us parametrize $\mu$ using the spherical angles: $\mu = \mu(\phi)$, where $\phi = (\phi_1, \ldots, \phi_{p-1})^t$.
Let $X_1, \ldots X_n$ be a random sample from the Watson distribution.

If we take the sample to be fixed and the underlying parameters to be unknown, then the log-likelihood function of \eqref{eqn:watson} is given by:
\begin{equation*}
  \ell(\phi,\kappa) = n \log d(\kappa) + \kappa \sum_{j=1}^n (X_j^t \mu(\phi))^2.
\end{equation*}
The maximum likelihood estimation of $\mu$ and $\kappa$ comes from the estimating equation:
\begin{equation} \label{eqn:gradient}
  \nabla_{\phi,\kappa} \ell(\phi,\kappa) = 0,
\end{equation}
where $\nabla_{\phi,\kappa}$ denotes the gradient.
Let $\hat{\phi}$ and $\hat{\kappa}$ be the solutions to \eqref{eqn:gradient}, which are the maximum likelihood estimators.
Then the standard theory of maximum likelihood estimators~\cite[pp.294-296]{coxHinkley:book} shows that the large sample asymptotics satisfy:
\begin{equation}
  \label{eqn:asymptotics}
  \sqrt{n}\left[ \left( \begin{array}{c} \hat{\phi} \\ \hat{\kappa} \end{array} \right) - \left( \begin{array}{c} {\phi} \\ {\kappa} \end{array} \right) \right] \to_d N_p(0,I(\phi,\kappa)^{-1})
\end{equation}
as $n \to \infty$, where ``$\to_d$'' means convergence in
distribution, $I(\phi,\kappa)$ is the Fisher information
matrix\footnote{The Fisher information matrix is defined to be
  $I(\phi,\kappa) = -E\nabla^2_{\phi,\kappa} \ell(\phi,\kappa)$, where
  $\nabla^2_{\phi,\kappa}$ is the $p \times p$ Hessian matrix.} and
$N_p$ stands for the $p$-dimensional normal distribution with given
mean and covariance.
It turns out that in the case of the Watson distribution,
\begin{equation*}
  I(\phi,\kappa) = \left[ \begin{array}{ccc} * & \vline & 0 \\ \hline 0 & \vline & -\frac{\partial^2}{\partial \kappa^2} \log d(\kappa) \end{array} \right].
\end{equation*}
Consequently, from~\eqref{eqn:asymptotics}, we have that 
\begin{equation*}
  \sqrt{n} (\hat{\kappa} -\kappa) \to_d N_1 \left( 0, - \left( \frac{\partial^2}{\partial \kappa^2} \log d(\kappa) \right)^{-1} \right),
\end{equation*}
as $n \to \infty$.

\bibliographystyle{halpha}
\bibliography{my}

\end{document}